\documentclass[12pt]{article}

\usepackage{amsmath,amssymb,amsthm}

\begin{document}

\newcommand{\End}{{\rm{End}\ts}}
\newcommand{\non}{\nonumber}
\newcommand{\wt}{\widetilde}
\newcommand{\wh}{\widehat}
\newcommand{\ot}{\otimes}
\newcommand{\la}{\lambda}
\newcommand{\al}{\alpha}
\newcommand{\be}{\beta}
\newcommand{\ga}{\gamma}
\newcommand{\de}{\delta^{}}
\newcommand{\om}{\omega^{}}
\newcommand{\ib}{\bar{\imath}}
\newcommand{\jb}{\bar{\jmath}}
\newcommand{\kb}{\bar{k}}
\newcommand{\lb}{\bar{l}}
\newcommand{\hra}{\hookrightarrow}
\newcommand{\ve}{\varepsilon}
\newcommand{\ts}{\,}
\newcommand{\tss}{\hspace{1pt}}
\newcommand{\U}{ {\rm U}}
\newcommand{\Y}{ {\rm Y}}
\newcommand{\C}{\mathbb{C}\tss}
\newcommand{\Z}{\mathbb{Z}}
\newcommand{\ZZ}{{\rm{Z}}}
\newcommand{\A}{\mathcal{A}}
\newcommand{\gl}{\mathfrak{gl}}
\newcommand{\oa}{\mathfrak{o}}
\newcommand{\spa}{\mathfrak{sp}}
\newcommand{\g}{\mathfrak{g}}
\newcommand{\ka}{\mathfrak{k}}
\newcommand{\p}{\mathfrak{p}}
\newcommand{\sll}{\mathfrak{sl}}
\newcommand{\agot}{\mathfrak{a}}
\newcommand{\qdet}{ {\rm qdet}\ts}
\newcommand{\sdet}{ {\rm sdet}\ts}
\newcommand{\sgn}{ {\rm sgn}\ts}
\newcommand{\Sym}{\mathfrak S}

\renewcommand{\theequation}{\arabic{section}.\arabic{equation}}

\newtheorem{thm}{Theorem}[section]
\newtheorem{lem}[thm]{Lemma}
\newtheorem{prop}[thm]{Proposition}
\newtheorem{cor}[thm]{Corollary}

\theoremstyle{definition}
\newtheorem{defin}[thm]{Definition}
\newtheorem{example}[thm]{Example}

\theoremstyle{remark}
\newtheorem{remark}[thm]{Remark}

\newcommand{\bth}{\begin{thm}}
\renewcommand{\eth}{\end{thm}}
\newcommand{\bpr}{\begin{prop}}
\newcommand{\epr}{\end{prop}}
\newcommand{\ble}{\begin{lem}}
\newcommand{\ele}{\end{lem}}
\newcommand{\bco}{\begin{cor}}
\newcommand{\eco}{\end{cor}}
\newcommand{\bde}{\begin{defin}}
\newcommand{\ede}{\end{defin}}
\newcommand{\bex}{\begin{example}}
\newcommand{\eex}{\end{example}}
\newcommand{\bre}{\begin{remark}}
\newcommand{\ere}{\end{remark}}

\newcommand{\bal}{\begin{aligned}}
\newcommand{\eal}{\end{aligned}}
\newcommand{\beq}{\begin{equation}}
\newcommand{\eeq}{\end{equation}}
\newcommand{\ben}{\begin{equation*}}
\newcommand{\een}{\end{equation*}}

\newcommand{\bpf}{\begin{proof}}
\newcommand{\epf}{\end{proof}}

\def\beql#1{\begin{equation}\label{#1}}

\title{\Large\bf  Quasideterminants and Casimir elements for
the general linear Lie superalgebra}
\author{Alexander Molev\quad and\quad Vladimir Retakh}

\date{} 

\maketitle

\begin{abstract}
We apply the techniques of quasideterminants
to construct new families of Casimir elements for the general linear
Lie superalgebra $\gl(m|n)$ whose images
under the Harish-Chandra isomorphism are respectively
the elementary, complete and power sums supersymmetric functions.
\end{abstract}

\vspace{43 mm}

\noindent
School of Mathematics and Statistics\newline
University of Sydney,
NSW 2006, Australia\newline
alexm@maths.usyd.edu.au

\vspace{7 mm}

\noindent
Department of Mathematics\newline
Rutgers University,
Piscataway, NJ 08854, USA\newline
vretakh@math.rutgers.edu

\newpage

\section{Introduction}\label{sec:int}
\setcounter{equation}{0}

Let $A$ be a square matrix over a ring. Its {\it quasideterminants\/}
are certain rational expressions in the entries of $A$.
The theory of quasideterminants originates from the papers by
Gelfand and Retakh~\cite{gr:dm, gr:tn} and since then a
number of applications of the theory has been found;
see \cite{ggrw:q} for an overview. In particular,
the techniques of quasideterminants is fundamental
in the theory of noncommutative symmetric functions developed
by Gelfand, Krob, Lascoux, Leclerc, Retakh and Thibon~\cite{gkllrt:ns}.
The symmetric functions associated with a matrix whose entries
are elements of a noncommutative ring
is one of the interesting specializations of the general theory.
When applied to the matrix $E$ formed by the generators of the
general linear Lie algebra $\gl(n)$ the theory produces
a new family of Casimir elements for $\gl(n)$ as well as
a distinguished set of generators of the Gelfand--Tsetlin
subalgebra of $\U(\gl(n))$; see \cite[Section~7.4]{gkllrt:ns}.
These results were extended to the orthogonal and symplectic
Lie algebras in \cite{m:ns} with the use of the twisted Yangians
and quantum determinants; see also a review paper \cite{m:ya}.

In this paper we use the techniques of quasideterminants
to get new families of Casimir elements for the general linear
Lie superalgebra $\gl(m|n)$ and calculate their images with respect to
the Harish-Chandra isomorphism. They can be regarded as
super-analogs of those constructed in \cite[Section~7.4]{gkllrt:ns}.
Three families of Casimir elements are given explicitly in terms
of some oriented graphs associated with $\gl(m|n)$. The
Harish-Chandra images turn out to be respectively the elementary, complete
and power sums supersymmetric functions.

The starting point for our construction is a result of Nazarov~\cite{n:qb}.
He produced a formal series $B(t)$ called quantum Berezinian with coefficients in
the center of the universal enveloping algebra $\U(\gl(m|n))$.
Our first result is
a quasideterminant factorization of $B(t)$ (Theorem~\ref{thm:decomp}).
We then use it to get graph presentations for the Casimir elements
(Theorem~\ref{thm:casim}).

Some other families of Casimir elements for $\gl(m|n)$
were constructed e.g. in \cite{m:fs}. This work is
a super-version of
the earlier constructions of \cite{naz2,o:qi} for $\gl(n)$
and it provides a linear basis of the center of $\U(\gl(m|n))$
formed by the so-called quantum immanants.

\medskip

We acknowledge the financial support of the Australian Research Council.
The second author would like to thank the School of Mathematics and Statistics,
the University of Sydney, for the warm hospitality during his visit. His 
work was also partially supported by NSA.

\section{Preliminaries}\label{sec:pre}
\setcounter{equation}{0}

Let $x=(x_1,\dots,x_m)$ and $y=(y_1,\dots,y_n)$ be two
families of variables.
A polynomial $P$ in $x$ and $y$ is called
supersymmetric if $P$ is
symmetric separately in $x$ and $y$
and satisfies the following cancellation property:
the result of setting $x_m=-y_n=z$ in $P$ is independent of $z$.
We denote by $\Lambda(m|n)$ the algebra of supersymmetric
polynomials in $x$ and $y$. The algebra $\Lambda(m|n)$
is generated by the polynomials
\beql{powsum}
p_k=x_1^k+\cdots+x_m^k+(-1)^{k-1}(y_1^k+\cdots+y_n^k),\qquad k\geq 1,
\end{equation}
called the {\it power sums supersymmetric functions\/}.
Two other families of generators of
$\Lambda(m|n)$ are comprised by the {\it elementary\/}
and {\it complete\/} supersymmetric functions defined
respectively by the formulas
\beql{ec}
\bal
e_k=\sum_{p+q=k}
\sum_{i_1<\cdots <i_p}\sum_{j_1\leq\cdots\leq j_q}
x_{i_1}\cdots x_{i_p}\ts
y_{j_1}\cdots y_{j_q},\\
h_k=\sum_{p+q=k}
\sum_{i_1\leq\cdots \leq i_p}\sum_{j_1<\cdots< j_q}
x_{i_1}\cdots x_{i_p}\ts
y_{j_1}\cdots y_{j_q};
\eal
\end{equation}
see \cite{s:ce}, \cite{s:cs}.

We shall denote by $E_{ij}$, $i,j=1,\dots,m+n$ the standard basis of the
Lie superalgebra $\gl(m|n)$. The $\Z_2$-grading on $\gl(m|n)$ is
defined by $E_{ij}\mapsto \ib+\jb$, where $\ib$
is an element of $\Z_2$ which equals
$0$ or $1$
depending on whether $i\leq m$ or $i>m$. The commutation relations
in this basis are given by
\beql{commrel}
[E_{ij},E_{kl}]=\delta_{kj}E_{il}-\delta_{il}E_{kj}
(-1)^{(\ib+\jb)(\kb+\lb)}.
\end{equation}

Given a $m+n$-tuple
$(\la|\mu)=(\la_1,\dots,\la_m,\mu_1,\dots,\mu_n)\in\C^{m+n}$ we consider a
highest weight $\gl(m|n)$-module $L(\la|\mu)$ with the highest
weight $(\la|\mu)$. That is, $L(\la|\mu)$ is generated by a nonzero vector
$\xi$ such that
\beql{hwmod}
\bal
E_{ii}\ts\xi&=\la_i\ts\xi\qquad&&\text{for}\quad i=1,\dots,m,\\
E_{m+j,m+j}\ts\xi&=\mu_j\ts\xi\qquad&&\text{for}\quad j=1,\dots,n,\\
E_{ij}\ts\xi&=0\qquad&&\text{for}\quad 1\leq i<j\leq m+n.
\eal
\end{equation}
Any element $z$ of the center $\ZZ(\gl(m|n))$ of the universal
enveloping algebra $\U(\gl(m|n))$ acts in $L(\la|\mu)$ as a scalar
$\chi(z)$. For a fixed $z$ the scalar $\chi(z)$ is a
polynomial in $\la_i$ and $\mu_i$ which is supersymmetric
in the shifted variables defined by
\beql{shift}
\bal
x_i&=\la_{i}-i+1\qquad&&\text{for}\quad i=1,\dots,m,\\
y_j&=\mu_j+m-j\qquad&&\text{for}\quad j=1,\dots,n.
\eal
\end{equation}
Furthermore, the map $z\mapsto
\chi(z)$ defines an algebra isomorphism
\beql{hch}
\chi:\ZZ(\gl(m|n))\to \Lambda(m|n),
\end{equation}
which is called the Harish-Chandra isomorphism; see
\cite{k:rc}, \cite{s:ce}, \cite{s:cs}.

\section{Decomposition of the Quantum Berezinian}\label{sec:qber}
\setcounter{equation}{0}

Introduce the super-matrix $\wh E$ of size $(m+n)\times(m+n)$
whose $ij$-th entry is $\wh E_{ij}=(-1)^{\jb}_{}E_{ij}$.
By the {\it quantum Berezinian\/} we mean the formal series $B(t)$
defined by
\beql{qber}
\bal
B(t)=\sum_{\sigma\in S_m}\sgn\sigma{\ts}{}&
\big(1+{t}\ts{\wh E}\big)^{}_{\sigma(1),1}\cdots
\big(1+t\ts(\wh E-m+1)\big)^{}_{\sigma(m),m}\\
\times \sum_{\tau\in S_n}\sgn\tau{\ts}{}&
\big(1+t\ts(\wh E-m+1)\big)^{-1}_{m+1,m+\tau(1)}\cdots
\big(1+t\ts(\wh E-m+n)\big)^{-1}_{m+n,m+\tau(n)}.
\eal
\end{equation}
The quantum Berezinian was constructed by Nazarov~\cite{n:qb}.
He also proved that all its coefficients are central
in the universal enveloping algebra $\U(\gl(m|n))$.
The image of $B(t)$ under the Harish-Chandra isomorphism
is given by
\beql{berhch}
\chi\big(B(t)\big)=\frac{(1+tx_1)\cdots (1+tx_m)}{(1-ty_1)\cdots (1-ty_n)},
\end{equation}
cf. \cite{m:fs}. Our first result is a decomposition of $B(t)$ into a product of
quasideterminants.
If $X$ is a square matrix over a ring with $1$ such that there exists
the inverse matrix $X^{-1}$ and its $ji$-th entry $(X^{-1})_{ji}$
is an invertible element of the ring, then the $ij$-{\it th quasideterminant
of} $X$ is defined by the formula
\ben
|X|_{ij}=   \big((X^{-1})_{ji}\big)^{-1},
\end{equation*}
see \cite{gr:dm, gr:tn} for other  equivalent definitions
of the quasideterminants and
their properties.

\bth\label{thm:decomp}
We have the following decomposition of $B(t)$ in the algebra
of formal series with coefficients in $\U(\gl(m|n))$
\begin{multline}\label{decomp}
B(t)=\big|1+t\ts\wh E^{(1)}\big|_{11}\cdots
\big|1+t\ts(\wh E^{(m)}-m+1)\big|_{mm}\\[1em]
{}\times\big|1+t\ts(\wh E^{(m+1)}-m+1)\big|^{-1}_{m+1,m+1}\cdots
\big|1+t\ts(\wh E^{(m+n)}-m+n)\big|^{-1}_{m+n,m+n},
\end{multline}
where ${\wh E}^{(k)}$ denotes the submatrix of $\wh E$ corresponding to
the first $k$ rows and columns. Moreover, the factors in the decomposition
are pairwise permutable.
\eth

\bpf
We employ a quasideterminant decomposition
of the quantum determinant for the Yangian $\Y(\gl(r))$.
The latter is the
associative algebra with the
generators $t_{ij}^{(1)},t_{ij}^{(2)},\dots$ where $1\leq i,j\leq r$
and the following defining relations
\beql{defrely}
[t_{ij}(u), t_{kl}(v)]=\frac{1}{u-v}\big(t_{kj}(u)t_{il}(v)-
t_{kj}(v)t_{il}(u)\big),
\end{equation}
where
\beql{tij}
t_{ij}(u)= \delta_{ij} + t^{(1)}_{ij} u^{-1} + t^{(2)}_{ij}u^{-2} +
\cdots \in \Y(\gl(n))[[u^{-1}]].
\end{equation}
Consider the quantum
determinant of the matrix $T(u)=\big[t_{ij}(u)\big]$ defined by
the following equivalent formulas
\beql{qdet}
\bal
\qdet T(u)&=\sum_{\sigma\in S_r} \sgn\sigma\cdot t_{\sigma(1),1}(u)\cdots
t_{\sigma(r),r}(u-r+1)\\
{}&=\sum_{\sigma\in S_r} \sgn\sigma\cdot t_{1,\sigma(1)}(u-r+1)\cdots
t_{r,\sigma(r)}(u).
\eal
\end{equation}
It is well-known that the coefficients of this series
are algebraically independent generators
of the center of the algebra $\Y(\gl(r))$; see e.g. \cite{mno:yc}
for a proof.
For $1\leq k\leq n$ denote by $T^{(k)}(u)$
the submatrix of $T(u)$ corresponding the
first $k$ rows and columns. We have the following
quasideterminant decomposition
of $\qdet T(u)$ in the algebra $\Y(\gl(m))[[u^{-1}]]$
\beql{decom}
\qdet T(u)=|T^{(1)}(u)|_{11}\cdots |T^{(m)}(u-m+1)|_{mm},
\end{equation}
where the factors are pairwise permutable;
see \cite{m:ns} and also \cite{gr:dm}, \cite{kl:mi} for analogous
decompositions in the case of noncommutative determinants of different types.
Now we apply the algebra homomorphism $\Y(\gl(m))\to\U(\gl(m|n))$
given by
\beql{hom1}
T(u)\mapsto 1+\wh E^{(m)}u^{-1}
\end{equation}
to \eqref{decom}, set $u=t^{-1}$ and multiply
both sides by $(1-t)\cdots (1-(m-1)t)$. This
will represent the first determinant factor in \eqref{qber}
as a product of quasideterminants which comprise the first $m$
factors in \eqref{decomp}; cf. \cite{m:ns}.

Now consider the second factor in \eqref{qber}.
We shall use the subscript $(k)$ of a matrix
to indicate its submatrix obtained by removing the first
$k-1$ rows and columns.
Here we need another version of the decomposition \eqref{decom}
given by
\beql{decom2}
\qdet T(u)=|T_{(1)}(u-n+1)|_{11}\cdots |T_{(n)}(u)|_{nn}.
\end{equation}
Apply another homomorphism  $\Y(\gl(n))\to\U(\gl(m|n))$
defined by
\beql{hom2}
T(u)\mapsto \big[(1+\wh E\ts u^{-1})^{-1}\big]_{(m+1)},
\end{equation}
(see \cite{n:qb}) to both sides of \eqref{decom2} with $\qdet T(u)$
expanded by the second formula in \eqref{qdet}.
Now observe that by the Inversion Theorem for quasiminors
\cite{gr:dm, gr:tn}, we have for any $k\in\{1,\dots,n\}$
\beql{inver}
\Big|\big[(1+\wh E\ts (u-n+k)^{-1})^{-1}\big]_{(m+k)}\Big|_{m+k,m+k}
=\big|1+\wh E^{(m+k)}\ts (u-n+k)^{-1}\big|_{m+k,m+k}^{-1}.
\end{equation}
To complete the argument,
it remains to set $u=t^{-1}+n-m$ and divide both sides
of the relation by the product $(1+t(1-m))\cdots (1+t(n-m))$.

Finally, note that the product of the first $m+n-1$ factors in \eqref{decomp}
coincides with the quantum Berezinian for the subalgebra
$\gl(m|n-1)$ of $\gl(m|n)$. Therefore the last factor in \eqref{decomp}
is permutable with the elements of $\gl(m|n-1)$ by the centrality
of the quantum Berezinian. The proof is completed by an obvious induction.
\epf

\section{Casimir elements}\label{sec:ce}
\setcounter{equation}{0}

Let $A=(A_{ij})$ be a square matrix of size $l\times l$
with entries from an arbitrary ring and
let $t$ be a formal variable.
Fix an
integer $i$ between $1$ and $l$. Following \cite[Definition~7.19]{gkllrt:ns}
introduce the noncommutative symmetric functions
associated with the matrix $A$ and the index $i$ as follows.
The {\it elementary symmetric functions\/} $\Lambda_k^{(i)}$,
the {\it complete symmetric functions\/} $S_k^{(i)}$,
the {\it power sums symmetric functions of the first kind\/}
$\Psi_k^{(i)}$ and
the {\it power sums symmetric functions of the second kind\/}
$\Phi_k^{(i)}$ are defined by the formulas
\beql{esf}
\bal
1+\sum_{k=1}^{\infty}\Lambda_k^{(i)}\ts t^k&=|1+tA|_{ii},\\
1+\sum_{k=1}^{\infty}S_k^{(i)}\ts t^k&=|1-tA|_{ii}^{-1},\\
\sum_{k=1}^{\infty}\Psi_k^{(i)}\ts t^{k-1}&=
|1-tA|_{ii}\ts\frac{d}{dt}\ts|1-tA|_{ii}^{-1},\\
\sum_{k=1}^{\infty}\Phi_k^{(i)}\ts t^{k-1}&=
-\frac{d}{dt}\log\big(|1-tA|_{ii}\big).
\eal
\end{equation}
These functions are polynomials in the entries of
the matrix $A$ and can be interpreted in terms of graphs in the
following way. Let us consider the complete oriented graph $\A$ with
$l$ vertices $\{1,2,\dots,l\}$, the arrow from $i$ to $j$ being
labelled by $A_{ij}$. Then every path in the graph going from $i$
to $j$ defines a monomial of the form
$A_{ir_1}A_{r_1r_2}\cdots A_{r_{k-1}j}$. A {\it simple path} is
a path such that $r_s\ne i,j$ for every $s$. Then by
\cite[Proposition 7.20]{gkllrt:ns},
$(-1)^{k-1}\Lambda_k^{(i)}$ is the sum of all monomials labelling
simple paths in $\A$ of length $k$ going from $i$ to $i$;
$S_k^{(i)}$ is the sum of all monomials labelling paths in
$\A$ of length $k$ going from $i$ to $i$;
$\Psi_k^{(i)}$ is the sum of all monomials labelling paths in
$\A$ of length $k$ going from $i$ to $i$, where the coefficient of
each monomial is the length of the first return to $i$;
$\Phi_k^{(i)}$ is the sum of all monomials labelling paths in
$\A$ of length $k$ going from $i$ to $i$, where the coefficient of
each monomial is the ratio of $k$ to the number of returns to $i$.

For any $i=1,\dots,m$ consider the matrix
$\wh E^{(i)}-i+1$ and the noncommutative
symmetric functions associated with this matrix
and the index $i$. We keep the above notation for these
functions.
Similarly, for any $j=1,\dots,n$
consider the matrix
$-\wh E^{(m+j)}+m-j$ and the noncommutative
symmetric functions associated with this matrix
and the index $m+j$. Again, we denote the functions by the same symbols
and distinguish them by the upper index $m+j$.

\bth\label{thm:casim}
The algebra $\ZZ(\gl(m|n))$ is generated by each of the families
\beql{casim}
\bal
\Lambda_k=&\sum_{i_1+\cdots+i_{m+n}=k}\Lambda_{i_1}^{(1)}\cdots
\Lambda_{i_m}^{(m)}\ts S_{i_{m+1}}^{(m+1)}\cdots
S_{i_{m+n}}^{(m+n)},\\
S_k=&\sum_{i_1+\cdots+i_{m+n}=k}S_{i_1}^{(1)}\cdots
S_{i_m}^{(m)}\ts\Lambda_{i_{m+1}}^{(m+1)}\cdots
\Lambda_{i_{m+n}}^{(m+n)},\\
\Psi_k=&\sum_{i=1}^m\Psi_k^{(i)}+(-1)^{k-1}\sum_{j=1}^n\Psi_k^{(m+j)},\\
\Phi_k=&\sum_{i=1}^m\Phi_k^{(i)}+(-1)^{k-1}\sum_{j=1}^n\Phi_k^{(m+j)},
\eal
\end{equation}
where $k=1,2,\dots$. Moreover, $\Psi_k=\Phi_k$ for any $k$, and
the Harish-Chandra images of these generators
are respectively
the elementary, complete and power sums supersymmetric functions,
\beql{hchecp}
\chi(\Lambda_k)=e_k,\qquad \chi(S_k)=h_k, \qquad \chi(\Psi_k)=p_k.
\end{equation}
\eth

\bpf
Introduce the generating functions for the supersymmetric polynomials
\eqref{powsum} and \eqref{ec} by
\beql{genf}
\bal
p(t)&=\sum_{k=1}^{\infty}p_k\ts t^{k-1},\\
e(t)&=1+\sum_{k=1}^{\infty}e_k\ts t^{k},\\
h(t)&=1+\sum_{k=1}^{\infty}h_k\ts t^{k}.
\eal
\end{equation}
These functions are related by
\beql{relf}
h(t)=e(-t)^{-1},\qquad p(t)=-\frac{d}{dt}\log e(-t)=e(-t)\frac{d}{dt}\ts e(-t)^{-1},
\end{equation}
see e.g. \cite{m:sf}. On the other hand, by Theorem~\ref{thm:decomp}
we have
\beql{eleml}
1+\sum_{k=1}^{\infty}\Lambda_k\ts t^{k}=B(t)
\end{equation}
which proves that the elements $\Lambda_k$ are central in $\U(\gl(m|n))$.
Moreover, $\chi(B(t))=e(t)$ due to \eqref{berhch}
and so $\chi(\Lambda_k)=e_k$.
The proof is completed by applying \eqref{relf} and taking into account
the fact that the factors in the decomposition \eqref{decomp}
are mutually permutable;
cf. the argument for the case of $\gl(n)$ \cite[Section~7.4]{gkllrt:ns}.
\epf

\bex
We have
\beql{phi1}
\bal
\Psi_1&=\sum_{i=1}^m(E_{ii}-i+1)+\sum_{j=1}^n(E_{m+j,m+j}+m-j),\\
\Psi_2&=\sum_{i=1}^m\Big((E_{ii}-i+1)^2+2\sum_{k=1}^{i-1} E_{ik}E_{ki}\Big)\\
{}&-\sum_{j=1}^n\Big((E_{m+j,m+j}+m-j)^2-2\sum_{l=1}^{m+j-1} (-1)^{\lb}
E_{m+j,l}E_{l,m+j}\Big).
\eal
\end{equation}
\eex


\begin{thebibliography}{99}

\bibitem{gkllrt:ns}
{I. M. Gelfand, D. Krob, A. Lascoux,
B. Leclerc, V. S. Retakh and J.-Y. Thibon},
{\it Noncommutative symmetric functions},
{Adv. Math.} {\bf 112} (1995), 218--348.

\bibitem{gr:dm} {I. M. Gelfand and V. S. Retakh}, {\it Determinants
of matrices over noncommutative rings}, {Funct. Anal. Appl.} {\bf 25}
(1991), 91-102.

\bibitem{gr:tn} {I. M. Gelfand and V. S. Retakh}, {\it A theory
of noncommutative determinants and characteristic functions of graphs},
{Funct. Anal. Appl.} {\bf 26} (1992), 1-20;  Publ. LACIM, UQAM, Montreal,
{\bf 14}, 1-26.

\bibitem{ggrw:q} {I. Gelfand, S. Gelfand, V. Retakh and R. Wilson},
{\it Quasideterminants},
preprint {\tt math.QA/0208146}.

\bibitem{k:rc}
{V. G. Kac}, {\it Representations of classical Lie
superalgebras}, {in:}
``Differential Geometry Methods in Mathematical Physics
II", (K. Bleuer, H. R. Petry, A. Reetz, Eds.),
Lecture Notes in Math., Vol. 676, pp. 597--626.
Springer-Verlag, Berlin/Heidelberg/New York, 1978.

\bibitem{kl:mi}
{D. Krob and B. Leclerc}, {\it Minor identities for
quasi-determinants and quantum determinants}, {Comm.
Math. Phys.} {\bf 169} (1995), 1--23.

\bibitem{m:sf}
{I. G. Macdonald},
{\it Symmetric functions and Hall polynomials},
Oxford University Press, 2nd edition
1995.

\bibitem{m:ns}
{A. I. Molev,}
{\it Noncommutative symmetric functions and Laplace operators
for classical Lie algebras}, {Lett. Math. Phys.}
{\bf 35} (1995), 135-143.

\bibitem{m:fs}
{A. I. Molev,}
{\it Factorial supersymmetric Schur functions and super Capelli identities},
in: ``Kirillov's Seminar on Representation Theory",
(G.~I.~Olshanski, Ed.), Amer. Math. Soc. Transl.,
Ser. 2, Vol. 181, AMS, Providence, R.I., 1998, pp. 109--137.

\bibitem{m:ya}
{A. I. Molev},
{\it Yangians and their applications},
in: ``Handbook of Algebra", Vol. 3, (M. Hazewinkel, Ed.), Elsevier, 2003.

\bibitem{mno:yc}
{A. Molev, M. Nazarov and G. Olshanski},
{\it Yangians and classical Lie algebras},
Russian Math. Surveys
{\bf 51}
(1996),
205--282.

\bibitem{n:qb}
{M. L. Nazarov},
{\it Quantum Berezinian and the classical Capelli identity},
{Lett. Math. Phys.}
{\bf 21}
(1991),
123--131.

\bibitem{naz2}
{M. L. Nazarov},
{\it Yangians and Capelli identities},
in: ``Kirillov's Seminar on Representation Theory",
(G.~I.~Olshanski, Ed.), Amer. Math. Soc. Transl.,
Ser. 2, Vol. 181, AMS, Providence, R.I., 1998, pp. 139--163.


\bibitem{o:qi}
{A. Okounkov},
{\it Quantum immanants and higher Capelli identities},
{Transformation Groups}
{\bf 1}
(1996),
99--126.

\bibitem {s:ce}
{M. Scheunert},
{\it Casimir elements of Lie superalgebras},
{in:} ``Differential Geometry Methods in Mathematical Physics",
pp. 115--124,
Reidel, Dordrecht, 1984.

\bibitem {s:cs} {J. R. Stembridge},
{\it A characterization of supersymmetric
polynomials}, {J. Algebra} {\bf 95} (1985), 439--444.




\end{thebibliography}
\end{document}